\documentclass{amsart}
\title{Planar subspaces are intrinsically CAT(0)}
\author{Russell Ricks}
\address{Binghamton University, Binghamton, New York, USA}
\email{rmricks@umich.edu}

\keywords{CAT(0), curvature, planar}
\subjclass[2010]{53C45, 53C22}

\usepackage{amssymb}
\def\R{\mathbb R}

\DeclareMathOperator{\p}{par}

\newcommand{\slashfrac}[2]{#1 / #2}

\usepackage{graphicx}
\usepackage{hyperref}

\begin{document}

\theoremstyle{plain}
\newtheorem{thm}{Theorem}[section]
\newtheorem{theorem}[thm]{Theorem}
\newtheorem{maintheorem}{Theorem}
\renewcommand{\themaintheorem}{\Alph{maintheorem}}
\newtheorem{mainlemma}[maintheorem]{Lemma}
\renewcommand{\themainlemma}{\Alph{mainlemma}}
\newtheorem*{theorem*}{Theorem}
\newtheorem{cor}[thm]{Corollary}
\newtheorem{corollary}[thm]{Corollary}
\newtheorem{lemma}[thm]{Lemma}
\newtheorem{conj}[thm]{Conjecture}
\newtheorem{prop}[thm]{Proposition}
\newtheorem*{fact*}{Fact}
\newtheorem*{question}{Question}

\theoremstyle{remark}
\newtheorem*{rem}{Remark}
\newtheorem*{note}{Note}

\theoremstyle{definition}
\newtheorem*{defn}{Definition}
\newtheorem*{convention}{Convention}
\newtheorem{alg}[thm]{Algorithm}

\thanks{This paper was originally written as part of a masters thesis at Brigham Young University in spring 2010, and posted on the arXiv that same year.
Recent work by Lytchak and Wenger \cite{lytchak-wenger18} convinced the author of continued interest in such problems, and that it was time to formally publish the result.
A few minor changes were made from the arXiv version.
The author would therefore like to thank his masters thesis advisor, Eric Swenson, for helpful discussions and advice.
The author would also like to thank the anonymous referee for helpful suggestions.}

\maketitle

\begin{abstract}
Let $M_\kappa^2$ be the complete, simply connected, Riemannian
$2$-manifold of constant curvature $\kappa \le 0$.
Let $E$ be a closed, simply connected subspace of $M_\kappa^2$ with
the property that every pair of points in $E$ is connected by a
rectifiable path in $E$.  We show that under the induced path metric,
$E$ is a complete CAT($\kappa$) space.  We also show that the natural
notions of angle coming from the intrinsic and extrinsic metrics coincide
for all simple geodesic triangles.
\end{abstract}

\section{Introduction}
Let $M_\kappa^2$ be the complete, simply connected, Riemannian
$2$-manifold of constant curvature $\kappa \le 0$.  We show the
following.

\begin{maintheorem}[Theorem \ref{CAT(kappa)}]\label{main CAT(kappa)}
Let $E$ be a closed, simply connected subspace of $M_\kappa^2$ with
the property that every pair of points in $E$ is connected by a
rectifiable path in $E$.  Then $E$, under the induced path metric, is a complete {\rm CAT($\kappa$)} space.
\end{maintheorem}

Our method of proof is based on some angle calculations.
However, there are multiple angles relevant to the situation.
One natural notion of angle for triangles in $E$ comes from the subspace metric; we call this the \emph{extrinsic angle} (see Section \ref{sec:angles}).
Another natural angle, the well-known Alexandrov angle, comes from the induced path metric.
Many statements in CAT($\kappa$) geometry use the Alexandrov angle.
In our situation, the extrinsic angle is relatively easy to understand, but it is not clear \emph{a priori} that the two angles should coincide.
We prove that they are, in fact, equal.

\begin{maintheorem}[Theorem \ref{angles}]\label{main angles}
For every geodesic triangle in $E$ (under the induced path metric),
the extrinsic and Alexandrov angles are equal.
\end{maintheorem}

\begin{figure}[h]
\includegraphics{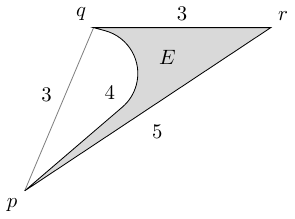}
\caption{The comparison angle in the induced path metric is $\angle^{(0)}_p(q,r) \approx 36.87^{\circ}$---strictly greater than the comparison angle in the subspace metric, which is $A_p(q,r) \approx 33.56^{\circ}$.}
\label{example}
\end{figure}

One might guess that for any given triangle, the comparison angle in the induced path metric is always smaller than or equal to the comparison angle in the subspace metric; then Theorem \ref{main CAT(kappa)} would follow immediately.
However, this guess is not correct, the following example shows:
Let $p,q,r \in \R^2$ ($= M_0^2$) have $d(p,r) = 5$ and $d(p,q) = d(q,r) = 3$.
Let $D \subset \R^2$ be the triangle with vertices $p,q,r$, along with its interior.
Let $L$ be a path in $D$ of length $4$ from $p$ to $q$.
Let $E \subset D$ be the closed disk with boundary $L$ and the line segments $[p,r]$ and $[q,r]$.
Then, for the geodesic triangle in $E$ with vertices $p,q,r$, the comparison angle in the induced path metric is $\angle^{(0)}_p(q,r) \approx 36.87^{\circ}$, whereas the comparison angle in the subspace metric is $A_p(q,r) \approx 33.56^{\circ}$ (see Figure \ref{example}).
Thus the comparison angle in the induced path metric can be strictly larger than the comparison angle in the subspace metric.

To further complicate the situation, one must prove the limit is correct for all possible approaches along the two geodesics (e.g.~ one cannot simply use the closest point on the other edge), and it is not obvious \emph{a priori} that the comparison angle in the induced path metric is decreasing along either edge.
(The inner comparison angle \emph{is} decreasing by the CAT(0) inequality, but that requires Theorem \ref{main CAT(kappa)}.)

We remark that a weaker version of Theorem \ref{main CAT(kappa)} is found elsewhere in the literature.
In fact, our proof is similar to
\cite{bishop} (where $E$ is the set of finite-distance points in the homeomorphic image of a closed disk), but the proof given there does not address the possibility that the extrinsic and Alexandrov angles may not coincide.

We also note that since this paper first appeared on the arXiv in spring 2010, a simple and elegant proof of Theorem \ref{main CAT(kappa)} has emerged:
In \cite[Proposition 12.1]{lytchak-wenger18},%
\footnote{stated for the case $\kappa = 0$}
Lythak and Wenger show that the metric completion $\bar{Y}$ of every locally-CAT($\kappa$) length space $Y$ (with $\kappa \le 0$) that is homeomorphic to the open disk, is CAT($\kappa$);
the idea of proof is simply to write $Y$ as an increasing limit of Jordan polygons in $Y$, each of which is easily seen to be locally CAT($\kappa$) and simply connected, hence $\bar{Y}$ is a limit of CAT($\kappa$) spaces and therefore CAT($\kappa$).
Now, under the hypotheses of Theorem \ref{main CAT(kappa)}, every simple geodesic triangle in $E$ has locally-CAT($\kappa$) interior and is homeomorphic to the open disk, hence the angles for simple geodesic triangles satisfy the CAT($\kappa$) inequality; the same statement readily follows for all geodesic triangles, and therefore $E$ is CAT($\kappa$) under the induced path metric.

Nevertheless, it is not clear how to obtain Theorem \ref{main angles} this way---neither directly from Theorem \ref{main CAT(kappa)}, nor using the method of proof sketched in the previous paragraph.
This paper provides a proof of Theorem \ref{main angles}.
And we use Theorem \ref{main angles} in a second paper \cite{ricks-complexes} to extend Theorem \ref{main CAT(kappa)} to the case of CAT($\kappa$) $2$-complexes.

As a final remark, one suspects that Theorems \ref{main CAT(kappa)} and \ref{main angles} remain valid when $\kappa > 0$, given appropriate diameter restrictions.
The proofs here could likely be adapted to this case, with suitable changes to deal with the loss of nonpositive curvature.
(One loses, for instance, convexity of the distance function; however, parallel lines can still be prevented from intersecting by restricting the diameter of the space.)

\section{Unique Geodesics}\label{sec:unique}
We make the following convention.

\begin{convention}
We will use the terms \emph{line} and \emph{line segment} to refer to standard geodesic lines and geodesic line segments in $M_\kappa^2$.  We will use \emph{geodesic} and \emph{geodesic segment} to refer to the geodesics and geodesic segments in $E$ under the induced path metric.
\end{convention}

Let $E$ be a closed, simply connected subspace of $M_\kappa^2$ with
the property that every pair of points in $E$ are connected by a
rectifiable path in $E$.  Let $d$ be the induced subspace metric and
$\bar d$ the induced path metric on $E$.  We will write $B_{d}(p,r)$
and $\overline B_{d}(p,r)$, respectively, for the open and closed
balls of radius $r$ about $p \in E$ in the standard metric on
$M_\kappa^2$.

Since $E$ is closed in $M_\kappa^2$, we know $(E,d)$ is complete.
The proofs of the following two more general results are provided
for completeness.

\begin{lemma}\label{complete}
The induced path metric on a complete metric space is complete.
\end{lemma}
\begin{proof}
Let $(X,d)$ be a complete metric space and $(X,\bar d)$ be the
induced path metric on $X$.  Suppose $\left\{ x_n \right\}_{n =
1}^{\infty}$ is a Cauchy sequence in $(X,\bar d)$.  Since $\bar
d(x,y) \ge d(x,y)$ for all $x,y \in X$, we know $\left\{ x_n
\right\}$ is also Cauchy in $(X,d)$.  Hence $x_n$ converges under
$d$ to some $x \in X$.  Now a Cauchy sequence converges if and only
if it has a convergent subsequence, so we may assume, by passing to
a subsequence if necessary, that $\bar d(x_n, x_m) < 2^{-m}$ for all
$m,n$ with $n > m$.  So for each $m$ there exists a path $c_m \colon
[0,1] \to X$ from $x_m$ to $x_{m+1}$ with $l(c_m) \le 2^{-m}$ by
assumption.  By linear reparameterization, we have paths $p_m \colon
[1-2^{-m+1},1-2^{-m}] \to X$ from $x_m$ to $x_{m+1}$ with $l(p_m)
\le 2^{-m}$.
Pasting these paths together and setting $p(1) = x$, we have a
continuous map $p \colon [0,1] \to X$.  Thus $p$ is a path from
$x_m$ to $x$ of length at most $\sum_{j=m}^{\infty} 2^{-j} =
2^{-m+1}$, so $\bar d(x_m, x) \le 2^{-m+1}$.  Therefore, $x_m \to x$
under $\bar d$.
\end{proof}

\begin{cor}
Suppose $X$ is a proper metric space, and every pair of points in $X$
is connected by a rectifiable path.  Then the induced path metric
on $X$ is geodesic.
\end{cor}
\begin{proof}
Let $x,y \in X$, and let $\delta$ be the distance from $x$ to $y$ in the
induced path metric.  Then by definition of path length, for each
$n \in \mathbb N$ there is a $(\delta + \slashfrac{1}{n})$-Lipschitz path
$\sigma_n \colon [0,1] \to X$ from $x$ to $y$.  By the Arzel\`a-Ascoli
theorem, these paths limit onto a $\delta$-Lipschitz path
$\sigma_0 \colon [0,1] \to X$ from $x$ to $y$.  This path realizes the shortest
length of a path from $x$ to $y$, hence it is geodesic.
\end{proof}

In Euclidean geometry, the following fact is often useful:  Given
any line $L$ and point $p$, there is a unique line $L'$ parallel to
$L$ that passes through the point $p$.  In hyperbolic geometry, we
no longer have a unique parallel line through $p$, so we choose a
nice one.

\begin{defn}  Let $L$ be a line and $p$ be a point in $M_\kappa^2$.
Let $K$ be the line segment from $p$ to the point $q \in L$ closest
to $p$.  There is a unique line $L'$ in $M_\kappa^2$ such that the
angle between $L'$ and $K$ is $\pi/2$.  We call $L'$ the \emph{closest
parallel to $L$ at $p$} and write $\p(L,p)$ for $L'$.
\end{defn}

\begin{lemma}\label{unique geodesics}
$(E,\bar d)$ is uniquely geodesic.
\end{lemma}
\begin{proof}
Suppose $\sigma \colon [a,b] \to E$ and $\tau \colon [a,b] \to E$
are distinct unit-speed geodesics with $p = \sigma(a) = \tau(a)$ and
$q = \sigma(b) = \tau(b)$.  Note that since both are unit-speed
geodesics, $\sigma(t)$ is in the image of $\tau$ if and only if
$\sigma(t) = \tau(t)$, and similarly for $\tau(t)$.  Since $\sigma$
and $\tau$ are distinct, there is some $t_0 \in (a,b)$ such that
$\sigma(t_0) \neq \tau(t_0)$, hence $\sigma(t_0)$ is not in the
image of $\tau$.  Taking the last $a' \in [a,t_0]$ and the first $b'
\in [t_0,b]$ such that $p' = \sigma(a')$ and $q' = \sigma(b')$ are
both in the image of $\tau$, we have that $C = \sigma([a',b']) \cup
\tau([a',b'])$ is a simple closed curve in $E$.
(See Figure \ref{fig:unique geodesics}.)

\begin{figure}[h]
\includegraphics[width=4.5in]{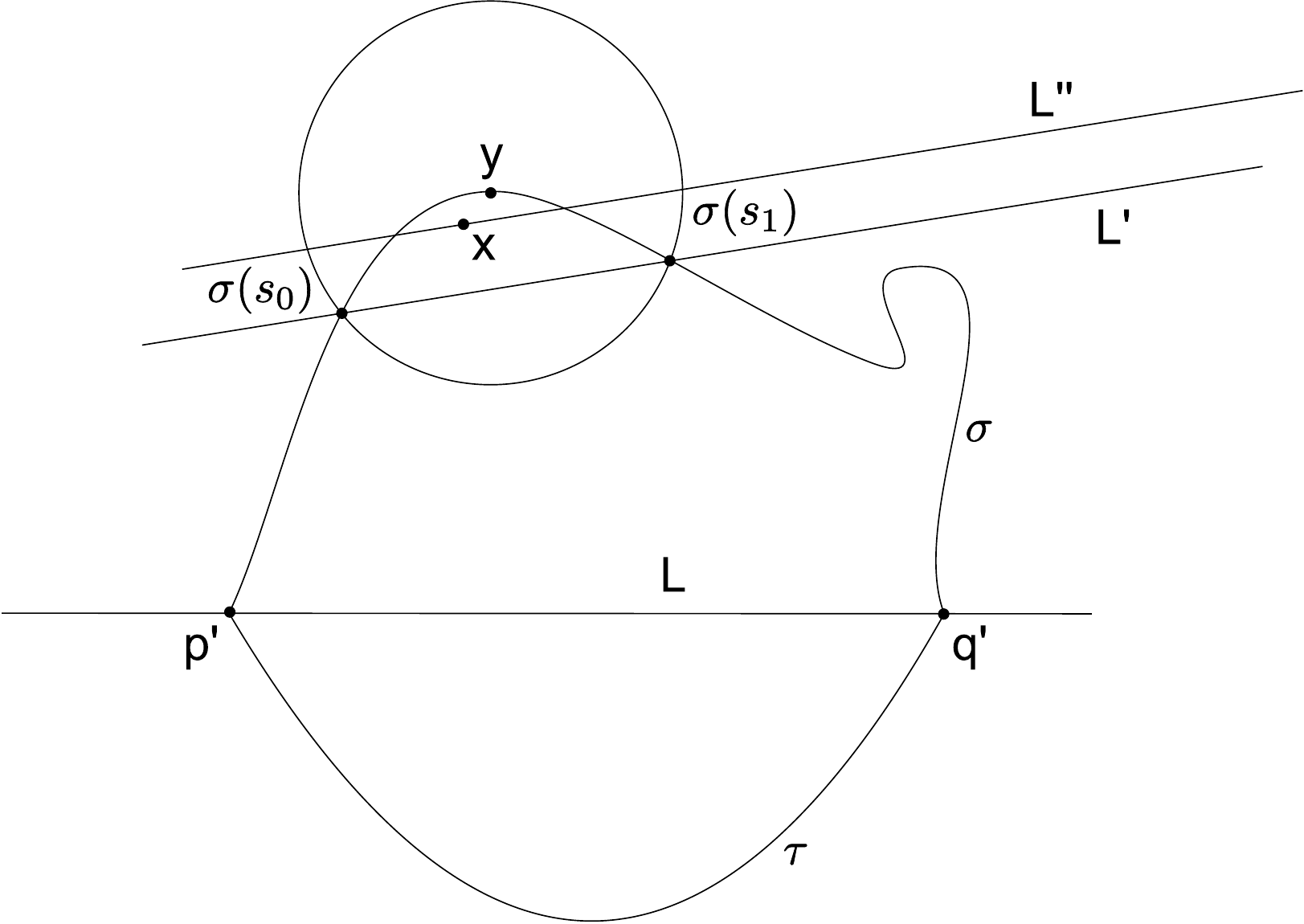}
\caption{Proving Lemma \ref{unique geodesics}.}
\label{fig:unique geodesics}
\end{figure}

Let $L$ be the line in $M_\kappa^2$ between $p'$ and $q'$.  Let $R$
be the maximum distance from $L$ to $C$, and let $t_1$ be the first
point of $[a',b']$ such that either $d(\sigma(t_1),L) = R$ or
$d(\tau(t_1),L) = R$.  We may assume $d(\sigma(t_1),L) = R$. Then,
since $C$ is a simple closed curve and $a' < t_1 < b'$, there is
some radius $r > 0$ about $y = \sigma(t_1)$ such that $\overline
B_d(y,r)$ does not intersect $\tau([a',b'])$.  Let $A$ be the
connected component of $C \cap \overline B_d(y,r)$ containing $y$,
and let $s_0 \in [a',t_1]$ and $s_1 \in [t_1,b']$ satisfy
$\sigma([s_0,s_1]) = A$.

Now let $L'$ be the line through $\sigma(s_0)$ and $\sigma(s_1)$.
Note that $d(\sigma(s_0),L) < d(y,L)$ and $d(\sigma(s_1),L) \le
d(y,L)$ by choice of $y$, so $y \notin L'$ by convexity of $d$.
By the Jordan curve theorem, $y$ is the limit of points in the
interior region $D$ bounded by $C$.  So there is some point $x \in
D$ with $d(x,y) < d(x,L')$.  Let $L'' = \p(L',x)$; since $d(x,L') =
d(L'',L')$, we also have $L' \cap L'' = \varnothing$.
Since $x$ is in $D$, $L''$ hits $C$ on each side of $x$; by
construction, $L''$ first hits $C$ inside $\overline B_d(y,r)$ in
each direction.  By choice of $r$, we therefore have a straight line
segment through $D$ between two points on $\sigma([a',b'])$ where
$\sigma$ does not follow the line segment exactly.  But $D \subset
E$ since $(E,d)$ is simply connected, so this contradicts $\sigma$
being geodesic.  Therefore, $(E,\bar d)$ is uniquely geodesic.
\end{proof}

\section{Simple Geodesic Triangles}\label{sec:triangles}
We will use the following terminology:  Call a geodesic triangle $T
\subset (E,\bar d)$ \emph{simple} if $T \subset (E, d)$ is a simple
closed curve.  For this section, let $T$
be a simple geodesic triangle in $(E,\bar d)$ with interior (under
the standard $M_\kappa^2$ metric) $S$ and exterior $U$.

\begin{prop}\label{convexity}
Let $L$ be a line in $M_\kappa^2$ that passes through two distinct
points $p$ and $q$ that lie on a single edge $A$ of $T$.  Let $L_0$
be the open line segment between $p$ and $q$.  If $L_0$ has empty
intersection with $T$ then $L_0 \subset U$.
\end{prop}
\begin{proof}
Since $T$ is a simple closed curve in $(E,d)$ and $(E,d)$ is simply
connected, $S \subset E$.  Hence if $L_0$ has empty intersection
with $T$, we have that $L_0$ is contained entirely in either $S$ or
$U$.  But $L_0 \subset S$ would give us $L_0 \subset E$, and this
contradicts the hypothesis that $A$ is the shortest path in $E$ from
$p$ to $q$.  Therefore, $L_0 \subset U$.
\end{proof}

\begin{figure}[h]
\includegraphics[width=2in]{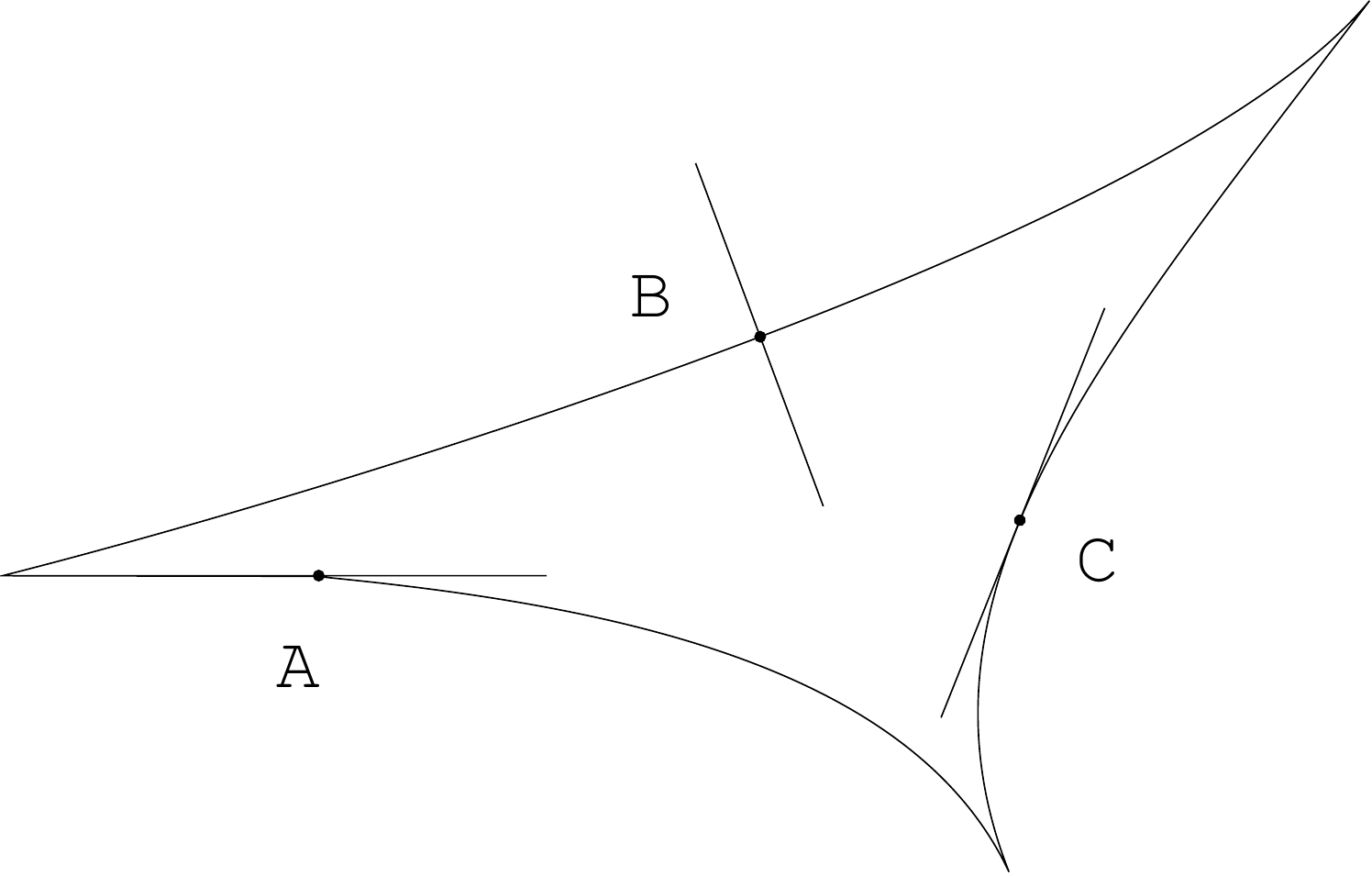}
\caption{Lemma \ref{interior/exterior} allows three types of lines
through an edge of a simple triangle:  (A) The line intersects the
triangle on one side and is locally interior or follows the edge on the other side, (B)
the line is locally interior on one side and locally exterior on the
other, or (C) the line is locally interior on both sides.}
\label{fig:interior exterior}
\end{figure}

\begin{lemma}\label{interior/exterior}
Let $L$ be a line in $M_\kappa^2$ that passes through the point $p
\in T$, where $p$ is not a vertex of $T$.  Let $A$ be the edge of
$T$ that contains $p$.  Suppose that $r > 0$ is a radius such that
$T \cap B_{d}(p,r) \subset A$, and let $L^-$ and $L^+$ be the two
components of $L \cap B_{d}(p,r) \setminus \left\{ p \right\}$.  Then
at least one of $L^-$ and $L^+$ has empty intersection with $U$.
Moreover, if $L^- \cap T
\neq \varnothing$ then $L^+
\cap U = \varnothing$.
(See Figure \ref{fig:interior exterior}.)
\end{lemma}

\begin{proof}
First suppose, by way of contradiction, that there exist points $x
\in L^- \cap U$ and $y \in L^+ \cap U$.  Let $r' > 0$ be some radius
with $r' < r$ such that we have both $B_{d}(x,r') \subset U$ and
$B_{d}(y,r') \subset U$.  Now by the Jordan Curve Theorem, $T =
\partial S$, so there is some point $q \in S$ close enough to $p$
that $L' = \p(L,q)$ hits points $x'$ in $B_{d}(x,r')$ and $y'$ in
$B_{d}(y,r')$.

Now $L'$ must be exterior at $x'$ and $y'$, but interior at $q$;
furthermore, $q$ lies between $x'$ and $y'$ on $L'$ by construction.
Thus $L'$ must hit $T$ somewhere between $x'$ and $q$ and somewhere
between $q$ and $y'$.  Therefore, $L'$ hits $T$ at two points $x''$
and $y''$ closest to $q$ (on opposite sides).  By hypothesis on the
radius $r$, we must have $x''\in A$ and $y'' \in A$.  Hence $L'$
contains a line segment between two points of $A$ that is completely
interior by construction.  This contradicts Proposition
\ref{convexity}, and therefore at least one of $L^-$ and $L^+$ has
empty intersection with $U$.

Suppose now that there is some point $z \in T \cap L^-$
and some point $w \in U \cap L^+$.  Let $r' > 0$ be some radius with
$r' < r$ such that we have $B_{d}(w,r') \subset U$.  The Jordan
Curve Theorem guarantees points in $U$ arbitrarily close to $z$, so
let $z' \in U$ be close enough to $z$ that the line $L''$ passing
through the points $z'$ and $p$ enters $B_{d}(w,r')$.  But then
$L''$ passes through the point $p$ and has nonempty intersection
with $U$ on both sides of $p$, which contradicts the result of the
previous paragraph.  Hence $L^+$ must have empty intersection with
$U$ if $L^-$ has nonempty intersection with $T$.
\end{proof}

\begin{cor}\label{three points}
Let $p_1$, $p_2$, and $p_3$ be three distinct points on a single
edge $A$ of $T$.  Suppose that $p_1$, $p_2$, and $p_3$ lie on a line
$L$ in $M_\kappa^2$, with $p_1$ and $p_3$ on opposite sides of
$p_2$.  Let $L_1$ and $L_2$ be the open line segments from $p_1$ to
$p_2$ and from $p_2$ to $p_3$, respectively.  If $L_1$ and $L_2$
both have empty intersection with $T \setminus A$, then the arc from
$p_1$ to $p_3$ along $T$ follows $L$.
\end{cor}
\begin{proof}
Suppose both $L_1$ and $L_2$ have empty intersection with $T
\setminus A$.  Then Proposition \ref{convexity} implies that both
$L_1$ and $L_2$ must have empty intersection with the interior.
Hence Lemma \ref{interior/exterior} gives us that if $L_1$ has
nonempty intersection with $U$,
then $L_2$ must follow $A$, so $L_2$ has nonempty intersection with
$T$, and thus $L_1$ has empty intersection with $U$; this is a
contradiction, so $L_1$ must have empty intersection with $U$.  Thus
$L_1$ follows $A$ (i.e., $L_1 \subset A$).  Similarly, $L_2$ must
follow $A$.  Therefore, the arc from $p_1$ to $p_3$ along $T$
follows $L$.
\end{proof}

\begin{lemma}\label{outer triangle}
Suppose the vertices of $T$ are $x$, $y$, and $z$.  Let $\triangle'$
be the triangle in $M_\kappa^2$ with vertices $x$, $y$, and $z$, and
let $C \subset M_\kappa^2$ be the convex hull of $\triangle'$.
Then $T$ is contained in $C$.
\end{lemma}

\begin{proof}
Suppose, by way of contradiction, that $p \in T \setminus C$.  Let
$L$ be the line passing through $x$ and $y$.  We may assume that $p$
lies in the component of $M_\kappa^2 \setminus L$ that contains no
point of $C$; let $H$ be the closure of this component.  Then $H
\cap T$ is compact and nonempty, so it contains at least one point
$p'$ of maximum distance to $L$.  Let $L'$ be the closest parallel to
$L$ at $p'$.  Now $L' \cap T$ is compact and nonempty, so let $q$ be
a point on $L' \cap T$ of maximum distance to $p'$.

Since $q \notin C$, $q$ is not a vertex of $T$.  Hence there is a
radius $r > 0$ such that $B_{d}(q,r)$ touches no point of any edge
of $T$ other than the one on which $q$ lies.  Let $L'^-$ and $L'^+$
be the two components of $L' \cap B_{d}(q,r) \setminus \left\{ q
\right\}$.  Lemma \ref{interior/exterior} requires both $L'^+$ and
$L'^-$ to be in $T$ since $L' \cap S$ is empty, but this contradicts
our choice of $q$.
Therefore, $T \subset C$, and the lemma is proved.
\end{proof}

\section{Extrinsic Angles}\label{sec:angles}
If $p$, $q$, and $r$ are distinct point in $E$, we will call the
angle in $M_\kappa^2$ at $p$ between $q$ and $r$ the \emph{outer
angle at $p$ between $q$ and $r$}, and denote it by $A_p (q,r)$.  Now
suppose $\sigma \colon [0,1] \to E$ and $\tau \colon [0,1] \to E$
are constant-speed geodesic line segments emanating from the point
$p \in E$, with $\sigma (1) = q$ and $\tau (1) = r$.  If the images
of $\sigma$ and $\tau$ intersect only at $p$, then by Proposition
\ref{convexity} and
Lemma \ref{outer triangle}, we have that $A_p (\sigma (t), \tau
(t'))$ decreases monotonically in both $t$ and $t'$, so the
\emph{extrinsic angle}
\[A'_p (q,r) = \lim_{t,t' \to 0} A_p (\sigma (t), \tau (t'))\]
is well defined.
If $\sigma \colon [0,1] \to E$ and $\tau \colon [0,1] \to E$
intersect at some $p' \neq p \in E$, then the limit above may not exist,
but $\sigma$ and $\tau$ must coincide from $p$ to $p'$ by Lemma
\ref{unique geodesics}; therefore we define
$A'_p (q,r) = 0$ in this case.

The concept of a CAT($\kappa$) space is closely related to the
Alexandrov angle at the vertex of a geodesic triangle.  Let
$\angle_p^{(\kappa)} (q,r)$ be the angle at $\bar p$ in the
comparison triangle $\triangle(\bar p, \bar q, \bar r)$ in
$M_\kappa^2$ for $\triangle(p, q, r)$.  The Alexandrov angle is
defined as
\[\angle_p (q,r) = \lim_{\epsilon \to 0} \sup_{0 < t,t' < \epsilon}
\angle_p^{(0)} (q,r).\]
We will show that the extrinsic angle $A'_p (q,r)$ equals the
Alexandrov angle $\angle_p (q,r)$.  We state the following two
results from
\cite{bridson}.  For more discussion on CAT($\kappa$) spaces, we
refer the reader to \cite{bridson} or \cite{ballmann}.

\begin{prop}[Paragraph preceding Definition I.2.15 of \cite{bridson}]
\label{alexandrov angle}
For any $\kappa \in \R$,
\[\angle_p (q,r) = \lim_{\epsilon \to 0} \sup_{0 < t,t' < \epsilon}
\angle_p^{(\kappa)} (q,r).\]
\end{prop}

\begin{prop}[Proposition I.1.14 of \cite{bridson}]
\label{triangle inequality}
Let $X$ be a metric space, and let $c,c',c'' \colon [0,1] \to X$ be geodesic segments issuing from the same point $p \in X$.  Then
\[\angle_p(c'(1),c''(1)) \le \angle_p(c(1),c'(1)) + \angle_p(c(1),c''(1)).\]
\end{prop}

To simplify the exposition, we again let $T$ be a simple geodesic triangle in $(E,\bar d)$
with interior (under the standard $M_\kappa^2$ metric) $S$ and
exterior $U$; denote the vertices by $p$, $q$, and $r$.  Also, let
$\sigma \colon [0,1] \to E$ and $\tau \colon [0,1] \to E$ be the
geodesic line segments from $p$ to $q$ and from $p$ to $r$,
respectively.

\begin{lemma}\label{support strut}
Suppose $A'_p (q,r) < \slashfrac{\pi}{2}$, and $\tau$ follows a line $L$
in $M_\kappa^2$ near $p$ (i.e., $\tau([0,\delta]) \subset L$ for
some $\delta > 0$).  Then there exists $t_1 > 0$ such that, for any
$t$ with $0 < t < t_1$, the line segment $L'$ from $\sigma(t)$ to the
closest point on $L$ lies in $S \cup T$.
\end{lemma}
\begin{proof}
Since $A_p (\sigma (t), \tau (t'))$ decreases monotonically in both
$t$ and $t'$, we may find some $\delta' \in (0,\delta]$ such that
$A_p (\sigma (t), \tau (t')) < \slashfrac{\pi}{2}$ for all $t$ and $t'$
with $0 < t, t' \le \delta'$.  Let $D = \overline
B_{d}(p,\epsilon)$, where $\epsilon > 0$ is small enough that $D
\cap T \subset \sigma([0,\delta']) \cup \tau([0,\delta'])$.  Let $P$
be the projection in $M_\kappa^2$ onto $L$, with domain restricted to
the image of $\sigma$, and let $L^+$ be the component of $L
\setminus \left\{ p \right\}$ that has nonempty intersection with
the image of $\tau$.

Since $A'_p (q,r) < \slashfrac{\pi}{2}$, there is some $t_0 > 0$ with $C
= \sigma([0,t_0]) \subset D$ such that $P(\sigma(t)) \in L^+$ for
every $t$ with $0 < t \le t_0$.  Since $P$ is continuous and $C$ is
compact, $P(C)$ has some point $q_1 = \sigma(t_1) \in C$ such that
$P(q_1)$ attains the maximum distance from $p$.  We further require
that $t_1$ be the smallest such value.
(See Figure \ref{fig:support strut}.)

\begin{figure}[h]
\includegraphics[width=4.5in]{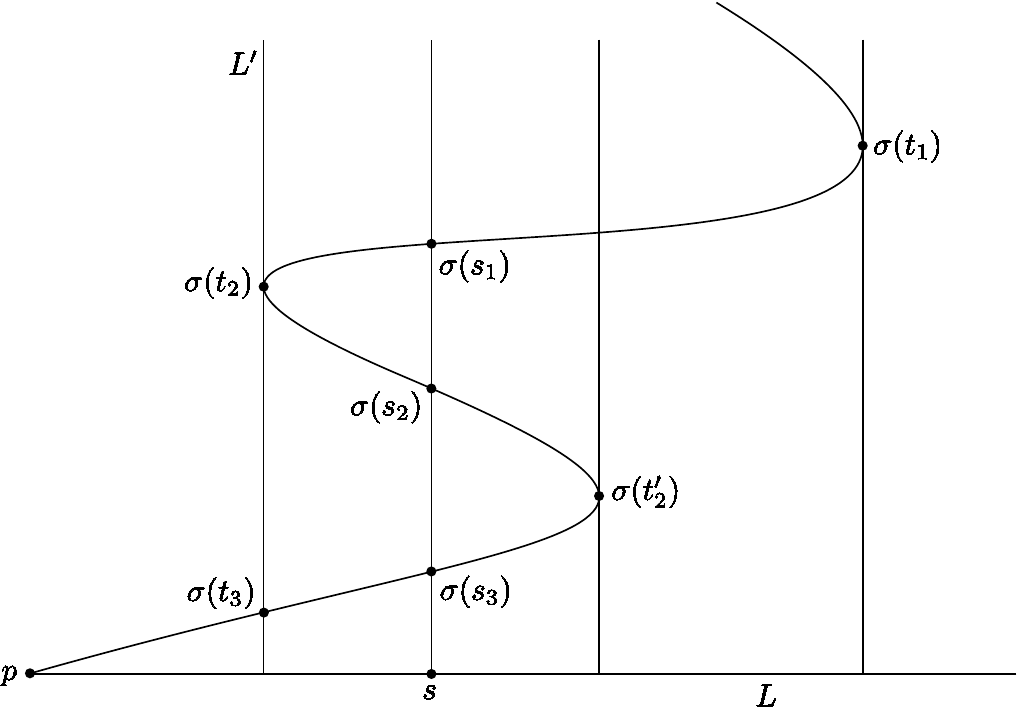}
\caption{Proving Lemma \ref{support strut}.}
\label{fig:support strut}
\end{figure}

Now suppose, by way of contradiction, the line segment $L'$ from
$q_2 = \sigma(t_2)$ to $P(q_2)$ contains a point of $U$ for some
$t_2$ with $0 < t_2 < t_1$ (note that $L' \perp L$).  Let $t_3$ be
the smallest positive value such that $q_3 = \sigma(t_3)$ lies on
$L'$.  If $t_3 = t_2$ then the line segment between $q_2$ and
$P(q_2)$ cuts one of $S$ or $U$ into two components; by Lemma
\ref{outer triangle}, it must therefore have interior in $S$, which
contradicts our hypothesis on $t_2$.  Thus $0 < t_3 < t_2$, and $L'$
has nontrivial intersection with $U$ between $q_2$ and $q_3$.  Hence
some $t'_2$ with $t_3 < t'_2 < t_2$ must have $P(\sigma(t'_2))$
farther from $p$ than $P(q_2) = P(q_3)$.  Let $q'_2 = \sigma(t'_2)$,
and let $s$ be the midpoint between $P(q_2)$ and $P(q'_2)$.  By the
intermediate value theorem, there must be some $s_1$ with $t_2 < s_1
< t_1$ such that $P(s_1) = s$.  Similarly, $P^{-1}(s)$ must contain
points $\sigma(s_2)$ and $\sigma(s_3)$ with $t'_2 < s_2 < t_2$ and
$t_3 < s_3 < t'_2$.  Thus these three points lie on a line in
$M_\kappa^2$ (orthogonal to $L$), so by Corollary \ref{three
points}, $q_2$ and $q'_2$ must also lie on this line; this is a
contradiction, so no such point $q_2$ can exist. Therefore, for any
$t$ with $0 < t < t_1$, the line segment from $\sigma(t)$ to $L$
perpendicular to $L$ is contained in $S \cup T$.
\end{proof}

\begin{lemma}\label{zero angles}
Suppose that $A'_p (q,r) = 0$ and $\tau$ follows a line $L$ in
$M_\kappa^2$ near $p$.  Then $\angle_p (q,r) = 0$.
\end{lemma}
\begin{proof}
For simplicity, we assume $\kappa = 0$ or $\kappa = -1$.  Let
$\epsilon > 0$ be given.  Since $A_p (\sigma (t), \tau (t'))$
decreases monotonically in both $t$ and $t'$, we may find some
$\delta > 0$ such that $A_p (\sigma (t), \tau (t')) < \epsilon$ for
all $t$ and $t'$ with $\sigma(t), \tau(t') \in \overline
B_{d}(p,\delta) \setminus \left\{ p \right\}$.
Write $D = \overline B_{d}(p,\delta)$.  Replacing $\delta$
by a smaller positive constant if necessary, we may assume that
every point of $T$ in $D$ is in the image of $\sigma$ or $\tau$ and
that the image of $\tau$ in $D$ follows $L$.  Let $P$ be the
projection from the image of $\sigma$ onto $L$, and let $t_1$ be the
point guaranteed by Lemma \ref{support strut}.

Let $\delta'$ be the distance in $M_\kappa^2$ from $p$ to
$P(\sigma(t_1))$, and note that $0 < \delta' < \delta$.  Suppose
that $q'$ and $r'$ are points in $B_{d}(p,\delta') \setminus \left\{
p \right\}$ along the images of $\sigma$ and $\tau$, respectively.
Let $a = d(p,q')$, $b = d(p,r')$, and $c = d(q',r')$, and let $\phi
= A_p (q', r')$.  Also let $a' = \bar d(p,q')$ and $c' = \bar
d(q',r')$; note that $a' \ge a$ and $c' \ge c$.  Since $\sigma$ is a
geodesic, the path straight from $p$ to $P(q')$ and then straight to
$q'$, which stays in $E$ by choice of $t_1$, must have length at
least $a'$.  Hence if $\kappa = 0$ then
\[a' \le a (\cos \phi + \sin \phi) \le a (1 + \sin \phi) \le
a (1 + \sin \epsilon) \le a (1 + \epsilon),\]
and if $\kappa = -1$ then by the hyperbolic law of sines,
\[\sinh a' \le (\cos \phi + \sin \phi) \sinh a
\le (1 + \epsilon) \sinh a.\]
Now suppose that $c' = c$.  By the law of cosines,
\[\cos \angle_p^{(0)} (q',r') = \frac{(a')^2 + b^2 -
c^2}{2(a')b} \ge \frac{a^2 + b^2 - c^2}{2(a')b} \ge \frac{a^2 + b^2
- c^2}{2a(1 + \epsilon)b} = \frac{1}{1 + \epsilon} \cos \phi,\]
and by the hyperbolic law of cosines,
\[\cos \angle_p^{(-1)} (q',r') =
\frac{\cosh a' \cosh b - \cosh c}{\sinh a' \sinh b} \ge \frac{\cosh
a \cosh b - \cosh c}{(1 + \epsilon) \sinh a \sinh b} =
\frac{1}{1 + \epsilon} \cos \phi.\]

On the other hand, suppose $c' > c$.
We claim $a > c$.
Let $L'_0$ be the line segment from $\sigma(t_1)$ to $P(\sigma(t_1))$, and let $L''$ be the line passing through both $p$ and $q'$.
By choice of $t_1$, the geodesic triangle $T_1$ with vertices $p$, $\sigma(t_1)$, and $P(\sigma(t_1))$ is simple.
Hence by Lemma \ref{outer triangle}, the segment of $L''$ that starts at $q'$ and ends on $L'_0$ must lie inside $S_1 \cup T_1$, where $S_1$ is the interior of $T_1$.
We also know the line segment from $q'$ to $P(q')$ lies in $S_1 \cup T_1$ by
choice of $t_1$.  Thus,
if $P(q')$ lies between $p$ and $r'$ on the line $L$, then the line
segment from $r'$ to $q'$ is contained in $S_1 \cup T_1$, which makes
$c' = c$.  But we are taking the case $c' > c$.  Therefore,
the outer angle $A_{r'} (p, q')$ is greater than $\slashfrac{\pi}{2}$,
and so $a > c$.

\begin{figure}[h]
\includegraphics{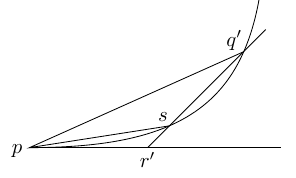}
\caption{Proving Lemma \ref{zero angles}.}
\label{fig:zero angles}
\end{figure}

Now consider the line segment in $M_\kappa^2$ from $r'$ to $q'$
(see Figure \ref{fig:support strut}):  It
hits $T_1$ at a first point $s$ (the edge hit is the one between $p$
and $q'$).  Let $\gamma$ be the path which travels from $r'$ to $s$
along the line segment and then from $s$ to $q'$ along $\sigma$.
Note that the length $\ell(\gamma)$ of $\gamma$ is at least $c'$.
Let $\alpha$ be the path that travels in a straight line from $p$ to
$s$ and then straight from $s$ to $q'$, and let $\alpha'$ be the
path that travels in a straight line from $p$ to $s$ and then from
$s$ to $q'$ along $\sigma$.  Note that $a \le \ell(\alpha) \le
\ell(\alpha') \le a'$.  Hence $a + c' \le \ell(\alpha) +
\ell(\gamma) = \ell(\alpha') + c \le a' + c$, and thus $a' - c' \ge
a - c$.  Therefore, $a > c$ gives us $a' - c' > 0$.  Since $a' \ge a
> 0$ and $c' \ge c > 0$, we have
\[(a')^{n + 1} - (c')^{n + 1} = (a' - c') \sum_{k = 0}^{n} (a')^k
(c')^{n - k} \ge
(a - c) \sum_{k = 0}^{n} a^k c^{n - k} = a^{n + 1} - c^{n + 1}\]
for all integers $n \ge 0$.  Hence for $\kappa = 0$ we have
\[\cos \angle_p^{(0)} (q',r') = \frac{(a')^2 + b^2 -
(c')^2}{2(a')b} \ge \frac{a^2 + b^2 - c^2}{2(a')b} \ge \frac{a^2 +
b^2 - c^2}{2a(1 + \epsilon)b} = \frac{1}{1 + \epsilon} \cos \phi,\]
and for $\kappa = -1$ we have
\begin{eqnarray*}
\cosh a' - \cosh c'
&=& \sum_{n = 1}^{\infty} \frac{1}{(2n)!} \left( (a')^{2n} - (c')^{2n}
\right) \\
&\ge& \sum_{n = 1}^{\infty} \frac{1}{(2n)!} \left( a^{2n} - c^{2n}
\right) \\
&=& \cosh a - \cosh c.
\end{eqnarray*}
Hence $\cosh a' - \cosh a \ge \cosh c' - \cosh c$, so the fact that
$\cosh b \ge 1$ gives us $\cosh a' \cosh b - \cosh a \cosh b \ge
\cosh c' - \cosh c$, and therefore $\cosh a' \cosh b - \cosh c' \ge
\cosh a \cosh b - \cosh c$.  Thus
\[\cos \angle_p^{(-1)} (q',r') =
\frac{\cosh a' \cosh b - \cosh c'}{\sinh a' \sinh b} \ge \frac{\cosh
a \cosh b - \cosh c}{(1 + \epsilon) \sinh a \sinh b} =
\frac{1}{1 + \epsilon} \cos \phi.\]
Thus, in either case,
\[\cos \angle_p^{(\kappa)} (q',r') \ge \frac{1}{1 + \epsilon}
\cos \phi,\]
and therefore we obtain $\angle_p (q',r') \le A'_p (q',r') = 0$ as
$\epsilon$ tends to zero.  This concludes the proof of the lemma.
\end{proof}

\begin{theorem}\label{angles}
In a simple geodesic triangle $T$ with vertices $p$, $q$, and $r$,
\[A'_p (q,r) = \angle_p (q,r).\]
That is, the extrinsic angle and the Alexandrov angle are equal.
\end{theorem}
\begin{proof}
By Proposition \ref{convexity}, the rays $R_{1,t}$ from $p$ through
$\sigma(t)$ limit monotonically to a ray $R_1$ as $t$ tends to zero.
Similarly, the rays $R_{2,t}$ from $p$ through $\tau(t)$
limit monotonically to a ray $R_2$ as $t$ tends to zero.

Suppose first that $R_1 \neq R_2$.  By construction, $R_1$ and $R_2$
are locally contained in $S \cup T$ near $p$.  Let $s_1$ be the last
point of $R_1$ contained in $S \cup T$.  Clearly, $s_1 \in T$; if
$s_1$ lies along $\sigma$ then $s_1$ must equal $q$ by Lemma
\ref{interior/exterior}.
Since $R_1$ is locally contained in $S \cup T$ near $p$, we have
$s_1 \neq p$, and thus $s_1$ cannot lie along $\tau$. Therefore,
$s_1$ lies along the geodesic segment between $q$ and $r$. Similarly,
the last point $s_2$ of $R_2$ that is contained in $S \cup T$ must
lie along the geodesic segment between $q$ and $r$.  Note that $\angle_p
(s_1,s_2) = A'_p (q,r)$, since both measure the angle between $R_1$
and $R_2$.

If $\sigma$ follows $R_1$ for some positive distance beyond $p$,
then $\angle_p (q,s_1) = 0$ by definition.  On the other hand, if
$\sigma$ does not follow $R_1$ for any positive distance beyond $p$,
then the geodesic triangle $T_1 = \triangle(p, q, s_1)$ is simple,
and $\angle_p (q,s_1) = 0$ by Lemma \ref{zero angles}.  Thus in
either case, $\angle_p (q,s_1) = 0$; similarly, $\angle_p (s_2,r) =
0$.  Hence
\[\angle_p (q,r) \le \angle_p (q,s_1) + \angle_p (s_1,s_2) +
\angle_p (s_2,r) = \angle_p (s_1,s_2)\]
and
\[\angle_p (s_1,s_2) \le \angle_p (s_1,q) +
\angle_p (q,r) + \angle_p (r,s_2) = \angle_p (q,r)\]
by Proposition \ref{triangle inequality}.  Therefore $\angle_p (q,r)
= \angle_p (s_1,s_2) = A'_p (q,r)$.

Finally, suppose $R_1 = R_2$; note that this gives $A'_p (q,r) = 0$.
If $\sigma$ follows $R_1$ for some positive distance beyond $p$,
then $\angle_p (q,r) = A'_p (q,r) = 0$ by Lemma \ref{zero angles}.
Thus we may assume, by symmetry, that neither $\sigma$ nor $\tau$
follows $R_1$ for any positive distance beyond $p$.  Then by
construction of $R_1 = R_2$, the last point $s$ of $R_1$ contained
in $S \cup T$ must be along the geodesic segment from $q$ to $r$.  Hence
the geodesic triangles $T_1 = \triangle(p, q, s)$ and $T_2 =
\triangle(p, s, r)$ are simple, and since $A'_p (q,s) = A'_p (s,r) =
0$ by construction, $\angle_p (q,s) = \angle_p (s,r) = 0$ by Lemma
\ref{zero angles}.  Therefore,
\[\angle_p (q,r) \le \angle_p (q,s) + \angle_p (s,r) = 0,\]
and the theorem is proved.
\end{proof}

\begin{theorem}\label{CAT(kappa)}
Let $M_\kappa^2$ be the complete, simply connected, Riemannian $2$-manifold of constant curvature $\kappa \le 0$.  Let $E$ be a closed, simply connected subspace of $M_\kappa^2$ with the property that every pair of points in $E$ is connected by a rectifiable path in $E$.  Then $E$, under the induced path metric $\bar d$, is a complete {\rm CAT($\kappa$)} space.
\end{theorem}
\begin{proof}
We want to show that $\angle_p (q,r) \le \angle_p^{(\kappa)} (q,r)$
for every triple of distinct points $p,q,r \in E$.  So let $T$ be
the geodesic triangle in $E$ with distinct vertices $p,q,r \in E$.
If the two geodesic segments from $p$ intersect at some point
$p' \neq p$, then
$\angle_p (q,r) = 0 \le \angle_p^{(\kappa)} (q,r)$.
On the other hand, if the two geodesic segments from $q$ intersect
at some point $q' \neq q$, then since the degenerate triangle with
vertices $p,q,q'$ automatically satisfies the CAT($\kappa$)
inequality, by Alexandrov's Lemma \cite[Lemma II.4.10]{bridson}
it suffices to prove the inequality for the case when these two
geodesic segments do not intersect except at $q$.
We may similarly assume the two geodesic segments from $r$ do not
intersect except at $r$.
Thus we may assume $T$ is simple.

As in the proof of Theorem \ref{angles}, we have two
(possibly equal) limit rays $R_1$ and $R_2$ from $p$.  Cutting along
these rays gives three (possibly degenerate) triangles.  The middle
triangle has Alexandrov angle at $p$ at most the extrinsic angle,
since only the edge opposite $p$ can be longer than the distance in
$M_\kappa^2$.  The two outside triangles have Alexandrov angle $0$
at $p$ by Lemma~ \ref{zero angles}.  So again by Alexandrov's Lemma \cite[Lemma II.4.10]{bridson}, we have $\angle_p (q,r) \le \angle_p^{(\kappa)} (q,r)$.
Therefore, $(E, \bar d)$ is CAT($\kappa$) \cite[Proposition II.1.7(4)]{bridson}.
(Completeness was proved in Lemma \ref{complete}.)
\end{proof}

\providecommand{\bysame}{\leavevmode\hbox to3em{\hrulefill}\thinspace}
\providecommand{\MR}{\relax\ifhmode\unskip\space\fi MR }
\providecommand{\MRhref}[2]{%
  \href{http://www.ams.org/mathscinet-getitem?mr=#1}{#2}
}
\providecommand{\href}[2]{#2}

\end{document}